\documentclass[11pt,twoside]{amsart}
\usepackage{latexsym,amssymb,amsmath}
\usepackage{curves}
\numberwithin{equation}{section}
\newtheorem{theorem}{Theorem}[section]
\newtheorem{lemma}[theorem]{Lemma}
\newtheorem{proposition}[theorem]{Proposition}
\newtheorem{corollary}[theorem]{Corollary}

\theoremstyle{definition}
\newtheorem{definition}[theorem]{Definition} 
 
\newtheorem{remark}[theorem]{Remark}
\newtheorem{example}[theorem]{Example}

\begin{document}


\title[On the vanishing ideal of an algebraic toric set]{On the
vanishing ideal of an algebraic toric set and its parameterized
linear codes}  

\author{Eliseo Sarmiento}
\address{
Departamento de
Matem\'aticas\\
Centro de Investigaci\'on y de Estudios
Avanzados del
IPN\\
Apartado Postal
14--740 \\
07000 Mexico City, D.F.
}
\email{esarmiento@math.cinvestav.mx}

\thanks{The first author was partially supported by CONACyT. The
second author is a member of the Center for Mathematical Analysis, 
Geometry and Dynamical Systems. 
The third author was partially supported by CONACyT 
grant 49251-F and SNI}

\author{Maria Vaz Pinto}
\address{
Departamento de Matem\'atica\\
Instituto Superior Tecnico\\
Universidade T\'ecnica de Lisboa\\ 
Avenida Rovisco Pais, 1\\ 
1049-001 Lisboa, Portugal 
}
\email{vazpinto@math.ist.utl.pt}

\author{Rafael H. Villarreal}
\address{
Departamento de
Matem\'aticas\\
Centro de Investigaci\'on y de Estudios
Avanzados del
IPN\\
Apartado Postal
14--740 \\
07000 Mexico City, D.F.
}
\email{vila@math.cinvestav.mx}


\begin{abstract} Let $K$ be a finite field and let $X$ be a subset of 
a projective space, over the field 
$K$, which is parameterized by monomials arising 
from the edges of a clutter.  We show some estimates for the 
degree-complexity, with respect
to the revlex order, 
of the
vanishing ideal $I(X)$ of $X$. If the clutter is uniform, we classify
the complete intersection property of $I(X)$ using linear algebra. We
show an upper bound for the minimum distance of 
certain parameterized linear codes along with certain estimates for
the 
algebraic invariants of $I(X)$.
\end{abstract}

\maketitle

\section{Introduction}

Let $K=\mathbb{F}_q$  be a finite field with $q\neq 2$ elements and 
let $y^{v_1},\ldots,y^{v_s}$ be a finite set of square-free monomials
with  
$s\geq 2$.  As usual if $v_i=(v_{i1},\ldots,v_{in})\in\mathbb{N}^n$, 
then we set 
$$
y^{v_i}=y_1^{v_{i1}}\cdots y_n^{v_{in}},\ \ \ \ i=1,\ldots,s,
$$
where $y_1,\ldots,y_n$ are the indeterminates of a ring of 
polynomials with coefficients in $K$. We shall always assume that
$\mathcal{A}=\{v_1,\ldots,v_s\}$ is the set of all characteristic
vectors of 
the edges of a clutter (see Definitions~\ref{clutter-def} and
\ref{charvec-def}). In
particular this means that the entries of $v_i$ are in $\{0,1\}$ for
all $i$. Consider the following set
parameterized  by these monomials 
$$
X:=\{[(x_1^{v_{11}}\cdots x_n^{v_{1n}},\ldots,x_1^{v_{s1}}\cdots
x_n^{v_{sn}})]\in\mathbb{P}^{s-1}	\vert\, x_i\in K^*\mbox{ for
all }i\}, 
$$
where $K^*=K\setminus\{0\}$ and $\mathbb{P}^{s-1}$ is a projective
space over the field $K$. The set $X$ is called an 
{\it algebraic toric set\/} parameterized  by
$y^{v_1},\ldots,y^{v_s}$. Let
$S=K[t_1,\ldots,t_s]=\oplus_{d=0}^\infty S_d$ 
be a polynomial ring 
over the field $K$ with the standard grading, let $[P_1],\ldots,[P_m]$
be the points of $X$, and 
let $f_0(t_1,\ldots,t_s)=t_1^d$. The {\it 
evaluation map\/} 
\begin{equation}\label{ev-map}
{\rm ev}_d\colon S_d=K[t_1,\ldots,t_s]_d\rightarrow K^{|X|},\ \ \ \ \ 
f\mapsto
\left(\frac{f(P_1)}{f_0(P_1)},\ldots,\frac{f(P_m)}{f_0(P_m)}\right) 
\end{equation}
defines a linear map of
$K$-vector spaces. The image of ${\rm ev}_d$, denoted by $C_X(d)$,
defines a {\it linear code}. Following \cite{algcodes} we call
$C_X(d)$ a {\it parameterized linear code\/} of
order $d$. As usual by a {\it linear code\/} we mean a linear
subspace of 
$K^{|X|}$. The {\it dimension\/} and the {\it length\/} of $C_X(d)$ 
are given by $\dim_K C_X(d)$ and $|X|$ respectively. The dimension
and the length 
are two of the {\it basic parameters} of a linear code. A third
basic parameter is the {\it minimum
distance\/} which is given by 
$$\delta_d=\min\{\|v\|
\colon 0\neq v\in C_X(d)\},$$ 
where $\|v\|$ is the number of non-zero
entries of $v$. The basic parameters of $C_X(d)$ are related by the
{\it Singleton bound\/} for the minimum distance:
$$
\delta_d\leq |X|-\dim_KC_X(d)+1.
$$

Parameterized linear codes are a nice subfamily of {\it evaluation
codes\/} (the notion of an evaluation code was introduced in
\cite{gold-little-schenck,GRT}). Parameterized linear codes were introduced and
studied in \cite{algcodes}. Some other families 
 of evaluation codes have been studied extensively
 \cite{delsarte-goethals-macwilliams,duursma-renteria-tapia,sorensen,tsfasman}.

The {\it vanishing ideal\/} of $X$, denoted by $I(X)$, is the
ideal of $S$ generated by the homogeneous polynomials of $S$ that
vanish on $X$. The contents of this paper are as follows. In
Section~\ref{prelim-invariants-of-I} we introduce the preliminaries 
and explain the connection between the invariants of the vanishing
ideal of $X$ and the parameters of $C_X(d)$. 

The ideal $I(X)$ is Cohen-Macaulay of height $s-1$
\cite{geramita-cayley-bacharach}. Recall that 
$I(X)$ 
is called a {\it complete
intersection\/} if it can be generated by $s-1$ homogeneous
polynomials of $S$. In \cite{ci-codes} it is shown that
$I(X)$ is a complete intersection if and only if $X$ is a projective
torus in $\mathbb{P}^{s-1}$ (see
Definition~\ref{projectivetorus-def}). If the clutter has all its 
edges of the same cardinality, in Section~\ref{ci-section} we
classify the complete  
intersection property of $I(X)$ using linear algebra (see
Theorem~\ref{ci-class-uniform}). 

Let $\succ$ be the reverse
lexicographical order on the monomials of $S$. Recall that the ideal
$I(X)$ has a 
unique reduced Gr\"obner basis with respect to $\succ$. The {\it
degree-complexity\/} of $I(X)$, with respect to $\succ$, is the
maximum degree in the reduced 
Gr\"obner basis of $I(X)$. In Section~\ref{d-compl-section} we study 
the structure of the reduced Gr\"obner basis of $I(X)$ 
and show an upper bound for the degree-complexity of $I(X)$ 
(see Theorem~\ref{gbasis-parameterized}). This means that the
algebraic methods of \cite{algcodes} to compute the invariants of
$I(X)$ will probably work better using the revlex order.  

In Section~\ref{upper-bounds-section} we show upper bounds for  
the minimum distance of $C_X(d)$ for a certain family of algebraic
toric sets $X$ arising from 
normal edge ideals (see Theorem~\ref{march14-10}(b)). For this family 
we also show estimates for the algebraic invariants of
$I(X)$.  The bounds on the minimum distance seem to indicate that 
the codes $C_X(d)$ that emerge from unicyclic connected graphs 
are especially nice from the point of view of their
error-correcting capacity and so are the codes $C_X(d)$ when $d$ is
small 
(see Remark~\ref{comparison-remark-1} and
Example~\ref{large-minimum-distance}).   
We give examples, within our family, of parameterized codes having a
large minimum distance 
relative to $|X|$ (see Example~\ref{large-minimum-distance}). Such
examples 
of linear codes with large minimum distance are essential, as they
show that 
our construction is attractive in the context of coding theory. The
codes $C_X(d)$ are only interesting when $d$ lies within a certain
range because $\delta_d=1$ for $d\gg 0$. This range is determined by 
${\rm reg}(S/I(X))$, the index of regularity of $S/I(X)$ 
(see Proposition~\ref{request-referee}). This is one of the
motivations to study the index of regularity. Another motivation 
comes from commutative algebra because, in our situation, ${\rm
reg}(S/I(X))$ is equal to the Castelnuovo-Mumford regularity which is 
an algebraic invariant of central importance in the area
\cite{eisenbud-syzygies}. The problem of finding a good decoding 
algorithm for our family of parameterized codes is not considered
here. The reader is referred to \cite[Chapter~9]{CLO1},
\cite{joyner-decoding,van-lint} and the references there for some
available 
decoding algorithms for some families of linear codes.

For all unexplained
terminology and additional information  we refer to
\cite{EisStu} (for the theory of binomial ideals),
\cite{CLO,Sta1} (for the theory of Gr\"obner bases and Hilbert
functions), and 
\cite{MacWilliams-Sloane,stichtenoth,tsfasman} (for the theory of 
error-correcting codes and algebraic geometric codes).

\section{Preliminaries}\label{prelim-invariants-of-I} 

We continue to use the notation and definitions used in the
introduction. In this section we introduce the basic algebraic
invariants of $S/I(X)$ and recall their connection with the basic 
parameters of parameterized linear codes. Then, we present a result on 
complete intersections that will be needed later.  

Recall that the {\it projective space\/} of 
dimension $s-1$ over $K$, denoted by 
$\mathbb{P}^{s-1}$, is the quotient space 
$$(K^{s}\setminus\{0\})/\sim $$
where two points $\alpha$, $\beta$ in $K^{s}\setminus\{0\}$ 
are equivalent if $\alpha=\lambda{\beta}$ for some $\lambda\in K$. We
denote the  
equivalence class of $\alpha$ by $[\alpha]$. Let
$X\subset\mathbb{P}^{s-1}$ be an algebraic toric set
parameterized by $y^{v_1},\ldots,y^{v_s}$ and let $C_X(d)$ be a
parameterized code of order $d$. The kernel of the
evaluation map ${\rm ev}_d$, defined in 
Eq.~(\ref{ev-map}), is precisely $I(X)_d$ the degree
$d$ piece of $I(X)$. Therefore there is an isomorphism of $K$-vector
spaces 
\begin{equation}\label{jan15-11}
S_d/I(X)_d\simeq C_X(d).
\end{equation}

Two of the basic parameters of $C_X(d)$ can be expressed 
using Hilbert functions of standard graded algebras \cite{algcodes,Sta1}, as we
explain below. Recall that the 
{\it Hilbert function\/} of
$S/I(X)$ is given by 
$$H_X(d):=\dim_K\, 
(S/I(X))_d=\dim_K\, 
S_d/I(X)_d.$$

The unique polynomial $h_X(t)=\sum_{i=0}^{k-1}c_it^i\in
\mathbb{Z}[t]$ of degree $k-1=\dim(S/I(X))-1$ such that
$h_X(d)=H_X(d)$ for 
$d\gg 0$ is called the {\it Hilbert polynomial\/} of $S/I(X)$. The
integer $c_{k-1}(k-1)!$, denoted by ${\rm deg}(S/I(X))$, is 
called the {\it degree\/} or  {\it multiplicity} of $S/I(X)$. In our
situation 
$h_X(t)$ is a non-zero constant because $S/I(X)$ has dimension $1$. 
Furthermore: 

\begin{proposition}{\rm(\cite[Lecture 13]{harris},
\cite{geramita-cayley-bacharach})} $h_X(d)=|X|$ for $d\geq |X|-1$.
\end{proposition} 

This result means that $|X|$ is equal to the {\it degree\/} 
of $S/I(X)$. From Eq.~(\ref{jan15-11}), we get the equality 
$H_X(d)=\dim_KC_X(d)$. Thus, we have: 

\begin{proposition}{\rm \cite{geramita-cayley-bacharach,GRT}}
$H_X(d)$ and ${\rm deg}(S/I(X))$ are equal to the
dimension and the length of $C_X(d)$ respectively. 
\end{proposition}

There are algebraic
methods, based on elimination theory and Gr\"obner bases, to compute
the dimension and the length of $C_X(d)$ \cite{algcodes}. 

The {\it index of regularity\/} of $S/I(X)$, denoted by 
${\rm reg}(S/I(X))$, is the least integer $p\geq 0$ such that
$h_X(d)=H_X(d)$ for $d\geq p$. The degree and the index of regularity
can be 
read off the Hilbert series as we now explain. The Hilbert series of
$S/I(X)$ can be 
written as
$$
F_X(t):=\sum_{d=0}^{\infty}H_X(d)t^d=\frac{h_0+h_1t+\cdots+h_rt^r}{1-t},
$$
where $h_0,\ldots,h_r$ are positive integers. Indeed
$h_i=\dim_K(S/(I(X),t_s))_i$ for $0\leq i\leq r$ and
$\dim_K(S/(I(X),t_s))_i=0$ for $i>r$.  
This follows from the fact that $I(X)$ is a Cohen-Macaulay lattice
ideal  of dimension $1$ \cite{algcodes} and by observing that $\{t_s\}$ is a regular
system of 
parameters for $S/I(X)$ (see \cite{Sta1}). The number 
$r$ is equal to the index of regularity of $S/I(X)$ and the degree of
$S/I(X)$ is equal to $h_0+\cdots+h_r$ (see \cite{Sta1} or 
\cite[Corollary~4.1.12]{monalg}).  

A good parameterized code should have large $|X|$ and  
with $\dim_KC_X(d)/|X|$ and $\delta_d/|X|$ as large as possible. 
The following result gives an indication of where to look for
non-trivial parameterized codes. Only the codes $C_X(d)$ with $1\leq
d<{\rm reg}(S/I(X))$ have the potential to be good linear codes. 
\begin{proposition}\label{request-referee} $\delta_d=1$ for
$d\geq{\rm reg}(S/I(X))$.  
\end{proposition}

\begin{proof} Since $H_X(d)$ is equal to the dimension of $C_X(d)$ and
$H_X(d)=|X|$ for $d\geq{\rm reg}(S/I(X))$, by a direct application of
the Singleton  
bound we get that $\delta_d=1$ for $d\geq{\rm reg}(S/I(X))$.
\end{proof}

The definition of $C_X(d)$ can be extended to any finite subset
$X\subset\mathbb{P}^{s-1}$ of a projective space over a field $K$
\cite{gold-little-schenck,GRT}. In this generality---the resulting
linear code---$C_X(d)$ is called an {\it
evaluation code\/} 
associated to $X$ \cite{gold-little-schenck}. It is also called  a 
{\it projective Reed-Muller code\/} over the set $X$ 
\cite{duursma-renteria-tapia,GRT}.  In this paper we will only
deal with parameterized codes over finite fields. 

The parameters of evaluation
codes associated to $X$ have 
been computed in a number of cases. If $X=\mathbb{P}^{s-1}$, 
the parameters of $C_X(d)$ are described in
\cite[Theorem~1]{sorensen}. If $X$ is the image of the affine space
$\mathbb{A}^{s-1}$ under the map $x\mapsto [(1,x)]$, the parameters 
of $C_X(d)$ are described in
\cite[Theorem~2.6.2]{delsarte-goethals-macwilliams}.  If $X$ is a
projective torus, the parameters of $C_X(d)$ are described in
\cite{duursma-renteria-tapia} and \cite{ci-codes}. In this paper we 
give upper bounds for the parameters of certain parameterized codes. 

As seen above, parameterized codes are a special type of evaluation
codes. 
What makes a parameterized 
code interesting is the fact that the vanishing ideal of $X$ is a
binomial ideal 
\cite{algcodes}, which allows the computation of the dimension and
length using the computer algebra system {\it Macaulay\/}$2$
\cite{mac2}. The  index of regularity of $S/I(X)$ can  also be computed
using  {\it Macaulay\/}$2$, which is useful to find genuine 
parameterized codes (see Proposition~\ref{request-referee}).

\begin{definition}\label{projectivetorus-def} The set 
$\mathbb{T}=\{[(x_1,\ldots,x_s)]\in\mathbb{P}^{s-1}\vert\, x_i\in
K^*\mbox{ for all }i\}$ is called a {\it projective torus} 
in $\mathbb{P}^{s-1}$. 
\end{definition}

An algebraic toric set is a multiplicative group under
componentwise multiplication. Thus, a projective torus is a
multiplicative group. For future reference we recall the following
result on complete 
intersections. 

\begin{proposition}{\rm\cite[Theorem~1,
Lemma~1]{GRH}}\label{ci-summary} If $\mathbb{T}$ is a projective torus 
in $\mathbb{P}^{s-1}$, then
\begin{itemize}
\item[(a)] $I(\mathbb{T})=(\{t_i^{q-1}-t_1^{q-1}\}_{i=2}^s)$. 
\item[(b)] $ F_\mathbb{T}(t)=(1-t^{q-1})^{s-1}/(1-t)^s$.  
\item[(c)] ${\rm reg}(S/I(\mathbb{T}))=(s-1)(q-2)$ and ${\rm
deg}(S/I(\mathbb{T}))=(q-1)^{s-1}$. 
\end{itemize}
\end{proposition}

\section{The complete intersection property of
$I(X)$}\label{ci-section}  

We continue to use the notation and definitions used in the
introduction and in the preliminaries. In this section, we use linear
algebra to give an structure theorem---valid for uniform clutters---for 
the complete intersection property of $I(X)$.  

\begin{definition}\label{clutter-def}
A {\it clutter\/} $\mathcal{C}$ is a family $E$ of subsets of a
finite ground set $Y=\{y_1,\ldots,y_n\}$ such that if $f_1, f_2 \in
E$, then $f_1\not\subset f_2$. The ground set $Y$ is called the {\em
vertex set} 
of $\mathcal{C}$ and $E$ 
is called the {\em edge set} of $\mathcal{C}$, they are denoted by
$V_\mathcal{C}$ 
and $E_\mathcal{C}$  respectively. 
\end{definition}

Clutters are special hypergraphs \cite{cornu-book} and are sometimes
called {\it Sperner families\/} in the literature. One 
important example of a clutter is a graph with the vertices and edges
defined in the  
usual way for graphs. 

\begin{definition}\label{charvec-def}
Let $\mathcal{C}$ be a clutter with vertex set
$V_\mathcal{C}=\{y_1,\ldots,y_n\}$ and let $f$ be an edge of
$\mathcal{C}$. The {\it characteristic vector\/} of $f$ is the vector
$v=\sum_{y_i\in f}e_i$, where
$e_i$ is the $i${\it th} unit vector in $\mathbb{R}^n$. 
\end{definition}

Throughout
this paper we assume that $\mathcal{A}:=\{v_1,\ldots,v_s\}$ is the
set of all 
characteristic vectors of the edges of a clutter $\mathcal{C}$.

\begin{definition} If $a\in {\mathbb R}^s$, its {\it support\/} is
defined as ${\rm supp}(a)=\{i\, |\, a_i\neq 0\}$.  
Note that $a=a^+-a^-$, 
where $a^+$ and $a^-$ are two non negative vectors 
with disjoint support called the {\it positive\/} and {\it
negative\/} part of $a$ respectively.  
\end{definition}

\begin{lemma}\label{jan6-10-1} 
Let $\mathcal{C}$ be a clutter and let $f\neq 0$ be a homogeneous
binomial 
of $I(X)$ of the form $t_i^b-t^c$ with $b\in\mathbb{N}$,
$c\in\mathbb{N}^s$ and $i\notin{\rm supp}(c)$. Then 
\begin{enumerate}
\item[\rm(a)] $\deg(f)\geq q-1$.
\item[\rm(b)] If $\deg(f)=q-1$, then $f=t_i^{q-1}-t_j^{q-1}$ for some
$j\neq i$. 
\end{enumerate}
\end{lemma}

\begin{proof} For simplicity of notation assume that
$f=t_1^b-t_2^{c_2}\cdots t_r^{c_r}$, where $c_j\geq 1$ for all $j$ 
and $b=c_2+\cdots+c_r$. Then
\begin{equation}\label{jan7-10-2}
(x_1^{v_{11}}\cdots x_n^{v_{1n}})^b=(x_1^{v_{21}}\cdots
x_n^{v_{2n}})^{c_2}\cdots (x_1^{v_{r1}}\cdots
x_n^{v_{rn}})^{c_r}\ \ \ \ \forall\, (x_1,\ldots,x_n)\in (K^*)^n,
\end{equation}
where $v_i=(v_{i1},\ldots,v_{in})$. Let $\beta$ be a generator of the
cyclic 
group $(K^*,\,\cdot\, )$.

(a) We proceed by contradiction. Assume
that $b<q-1$. First we claim that if $v_{1k}=1$ for some $1\leq k\leq
n$, then $v_{jk}=1$ for $j=2,\ldots,r$. To prove the claim assume that
$v_{1k}=1$ and $v_{jk}=0$ for some $j\geq 2$. Then, making $x_i=1$
for $i\neq k$ in 
Eq.~(\ref{jan7-10-2}), we get $(x_k^{v_{1k}})^b=x_k^b=x_k^m$, where 
$m=v_{2k}c_2+\cdots+v_{rk}c_r<b$. Then $x_k^{b-m}=1$ for $x_k\in
K^*$. In  particular  
$\beta^{b-m}=1$. Hence $b-m$ is a multiple of $q-1$ and consequently 
$b\geq q-1$, a contradiction. This completes the proof of the claim. 
Therefore ${\rm supp}(v_1)\subset{\rm supp}(v_j)$ for $j=2,\ldots,r$.
Since $\mathcal{C}$ is a clutter we get that $v_1=v_j$ for
$j=2,\ldots,r$, a contradiction because $v_1,\ldots,v_r$
are distinct. Hence $b\geq q-1$.

(b) It suffices to show that $r=2$. Assume $r\geq 3$.  We claim
that if $v_{2k}=1$ for some $1\leq k\leq 
n$, then $v_{jk}=1$ for $j\geq 3$. Otherwise, if $v_{2k}=1$ and
$v_{jk}=0$ for 
some $j\geq 3$, making $x_i=1$ for $i\neq k$ and $b=q-1$ in
Eq.~(\ref{jan7-10-2}) we get $1=x_k^m$ for any $x_k\in K^*$, for some
$0<m<q-1$. A contradiction because $\beta^m\neq 1$. 
This proves the claim. Therefore ${\rm supp}(v_2)\subset{\rm
supp}(v_j)$ for $j\geq 3$.   
As in part (a) we get $v_2=v_j$ for
$j\geq 3$, a contradiction. Hence $r=2$.
\end{proof}

The complete intersection property of $I(X)$ was first studied in
\cite{ci-codes}. We complement the following result by showing a
characterization of this property---valid for uniform clutters---using
linear algebra (see Theorem~\ref{ci-class-uniform}). 

\begin{theorem}{\cite{ci-codes}}\label{ci-characterization}
Let $\mathcal{C}$ be a clutter with $s$ edges and let 
$\mathbb{T}$ be a
projective torus in $\mathbb{P}^{s-1}$. The following are equivalent:
\begin{enumerate}
\item[($\mathrm{c}_1$)] $I(X)$ is a complete intersection.
\item[($\mathrm{c}_2$)]
$I(X)=(t_1^{q-1}-t_s^{q-1},\ldots,t_{s-1}^{q-1}-t_s^{q-1})$.
\item[($\mathrm{c}_3$)] $X=\mathbb{T}\subset\mathbb{P}^{s-1}$.
\end{enumerate}
\end{theorem}

For use below recall that the {\it toric ideal\/} associated to
$\mathcal{A}=\{v_1,\ldots,v_s\}$, denoted by $I_\mathcal{A}$, is the
prime ideal of 
$S=K[t_1,\ldots,t_s]$ given by (see \cite{Stur1}):
\begin{equation}\label{oct12-09}
I_\mathcal{A}=\left.\left(
t^a-t^b\right\vert\,
a=(a_i),b=(b_i)\in\mathbb{N}^s,\textstyle\sum_ia_iv_i=\sum_i
b_iv_i\right)\subset S.
\end{equation}

A clutter is called {\it uniform\/} if all its edges have 
the same number of elements. 

\begin{proposition}\label{ci+uniform-implies-li}
Let $\mathcal{C}$ be a uniform clutter. If $I(X)$ is a complete
intersection and $q\geq 3$, then $v_1,\ldots,v_s$ are linearly
independent. 
\end{proposition}

\begin{proof} To begin with we claim that if $f=t^{a^+}-t^{a^-}$ is
any 
non-zero homogeneous binomial in the lattice ideal $I(X)$, 
then 
$$
a=a^+-a^-\ \equiv\ 0\ {\rm mod} (q-1),
$$
that is, any entry of $a$ is a multiple of $q-1$. By
Theorem~\ref{ci-characterization} the degree of $f$ is at least $q-1$.
To show the claim we proceed by induction on $\deg(f)$. If 
$\deg(f)=q-1$, then by Theorem~\ref{ci-characterization} and
Lemma~\ref{jan6-10-1}(b) it is seen that $f=t_i^{q-1}-t_j^{q-1}$ for
some 
$i,j$, i.e., $a=(q-1)e_i-(q-1)e_j$. Assume that $\deg(f)>q-1$. By 
Theorem~\ref{ci-characterization} we obtain that $t^{a^+}$ and
$t^{a^-}$ are divisible by some $t_i^{q-1}$ and $t_j^{q-1}$
respectively. Then, $a_i^+\geq q-1$ and $a_j^-\geq q-1$ 
for some $i\in {\rm supp}(a^{+})$ and $j\in {\rm supp}(a^{-})$.
Therefore using that $f\in I(X)$ and the fact that $(K^*,\,\cdot\, )$
is a cyclic group 
of order $q-1$, it follows readily that the binomial 
\[
f'=\frac{t^{a^+}}{t_i^{q-1}}-\frac{t^{a^-}}{t_j^{q-1}}
\]   
is homogeneous, of degree $\deg(f)-(q-1)$, and belongs to $I(X)$.
Hence by induction hypothesis the vector
$(a^{+}-(q-1)e_i)-(a^{-}-(q-1)e_j)$ is a multiple 
of $q-1$, and so is $a=a^{+}-a^{-}$. This completes de proof of the
claim. 

To show that $v_1,\ldots,v_s$ are linearly independent we
proceed by contradiction. Assume that 
$v_1,\ldots,v_s$ are linearly dependent. As $\mathcal{C}$ is uniform,
there is a non-zero homogeneous binomial $f=t^{a^+}-t^{a^-}$ of least
degree in the toric ideal $I_\mathcal{A}$. This means that the degree
of $f$ is equal to 
the initial degree of $I_\mathcal{A}$ \cite[p.~ 110]{monalg}. Since
$I_\mathcal{A}\subset 
I(X)$ we obtain that $a=a^{+}-a^{-}$ is a multiple of $q-1$. 
Then, we can write $a^+=(q-1)b^{+}$, $a^-=(q-1)b^{-}$ for some
$b^+,b^-$ in $\mathbb{N}^s$. We set $u=t^{b^+}$, $v=t^{b^-}$, $g=u-v$,
$h=u^{q-2}+u^{q-3}v+\cdots+v^{q-2}$. From the equality 
$f=gh$ we obtain that $g\in I_\mathcal{A}$ or $h\in I_\mathcal{A}$
because $I_\mathcal{A}$ is a prime ideal and $q\geq 3$, a
contradiction to the choice 
of $f$ because $g$ and $h$ have degree less than that of $f$.
\end{proof}

\begin{definition} For  an ideal $I\subset S$ and  a polynomial $h\in
S$ the {\it 
saturation\/} of $I$ with respect to $h$
is the ideal
$$(I\colon
h^{\infty}):=\{f\in S\vert\, fh^m\in I\mbox{ for some }m\geq 1\}.
$$
We will only deal with the case where $h=t_1\cdots t_s$. 
\end{definition}

We call $\mathcal{A}$ {\it homogeneous\/} if $\mathcal{A}$ lies on an
affine 
hyperplane not containing the origin. Notice that if $\mathcal{C}$ is
uniform, then $\mathcal{A}$ is homogeneous. Given $\Gamma\subset
\mathbb{Z}^n$, the subgroup of $\mathbb{Z}^n$ generated by
$\Gamma$ will be denoted by $\mathbb{Z}\Gamma$. 

\begin{theorem}{\cite[Theorem~2.6]{algcodes}}\label{vila-dictaminadora}
Let $K=\mathbb{F}_q$ be a finite field, let 
$\mathcal{A}=\{v_1,\ldots,v_s\}\subset\mathbb{Z}^n$, and let
$\phi\colon
\mathbb{Z}^n/L\rightarrow\mathbb{Z}^n/L$
be the multiplication map $\phi(\overline{a})=(q-1)\overline{a}$,
where $L=\mathbb{Z}\{v_i-v_1\}_{i=2}^s$. If
$\mathcal{A}$ is homogeneous, then
\begin{equation}\label{happy-sept10-09}
((I_\mathcal{A}+(t_2^{q-1}-t_1^{q-1},\ldots,t_s^{q-1}-t_1^{q-1}))\colon
(t_1\cdots 
t_s)^\infty)\subset I(X)
\end{equation}
with equality if and only if the map $\phi$ is injective. 
\end{theorem}

We come to the main result of this section, a structure
theorem for complete intersections via linear algebra.

\begin{theorem}\label{ci-class-uniform} 
Let $\phi\colon
\mathbb{Z}^n/L\rightarrow\mathbb{Z}^n/L$
be the multiplication map $\phi(\overline{a})=(q-1)\overline{a}$,
where $L$ is the subgroup generated by $\{v_i-v_1\}_{i=2}^s$. If 
$\mathcal{C}$ is a uniform clutter and $q\geq 3$, then $I(X)$ is a
complete 
intersection if and only if $v_1,\ldots,v_s$ are linearly independent
and the map $\phi$ is injective.
\end{theorem}

\begin{proof} $\Rightarrow$) By
Proposition~\ref{ci+uniform-implies-li} the vectors 
$v_1,\ldots,v_s$ are linearly independent. Then $I_\mathcal{A}=(0)$
and 
by Theorem~\ref{ci-characterization} we
get the equality $I(X)=(\{t_1^{q-1}-t_i^{q-1}\}_{i=2}^s\})$. Hence, we
have equality in Eq.~(\ref{happy-sept10-09}). Therefore using
Theorem~\ref{vila-dictaminadora} we conclude that $\phi$ is 
injective. 

$\Leftarrow$) As the map $\phi$ is injective and $\mathcal{C}$ is
uniform, using Theorem~\ref{vila-dictaminadora},  we get the equality
$$((I_\mathcal{A}+(t_2^{q-1}-t_1^{q-1},\ldots,t_s^{q-1}-t_1^{q-1}))\colon
(t_1\cdots t_s)^\infty)=I(X).
$$
Since $\mathcal{A}$ is linearly 
independent one has that $I_\mathcal{A}=(0)$. Hence, the equality above
becomes $(\{t_1^{q-1}-t_i^{q-1}\}_{i=2}^s\})=I(X)$, i.e., $I(X)$ is a
complete intersection. 
\end{proof}

A graph with only one cycle is called {\it unicyclic}. 

\begin{corollary}\label{jan4-11} 
Let $\mathcal{C}$ be a unicyclic connected graph with $n$ vertices. 
If the only cycle of
$\mathcal{C}$ is odd, then $X=\mathbb{T}$ is a projective torus in
$\mathbb{P}^{n-1}$.  
\end{corollary}

\begin{proof} Assume that $\mathcal{C}$ is an 
odd cycle of length $n$. 
Let $y_1,\ldots,y_n$ be the vertices of $\mathcal{C}$. The
characteristic vectors of the edges of 
$\mathcal{C}$ are 
$$v_1=e_1+e_2,\, v_2=e_2+e_3,\ldots,v_{n-1}=e_{n-1}+e_n,\, v_n=e_n+e_1,$$
where $e_i$ is the $i${\it th} unit vector in $\mathbb{N}^n$. The
vectors $v_1,\ldots,v_n$ are linearly independent because $n$ is odd.
It is not hard to see that the quotient group
$\mathbb{Z}^n/\mathbb{Z}\{v_i-v_1\}_{i=2}^n$ is torsion-free. Hence,
by Theorem~\ref{ci-class-uniform}, $I(X)$ is a complete intersection.
Then, $X=\mathbb{T}$ is a projective torus in $\mathbb{P}^{n-1}$ by 
Theorem~\ref{ci-characterization}. If $\mathcal{C}$ is not an odd cycle,
then it has a vertex of degree $1$ and the proof follows by induction
because removing this vertex results in a graph that is connected and
has a unique odd cycle. 
\end{proof}

The next result shows that the index of regularity of complete
intersections 
associated to clutters provides an upper bound for the index of
regularity of $S/I(X)$. 

\begin{proposition}{\cite{ci-codes}}\label{reg-bound-loose} 
${\rm reg}(S/I(X))\leq (q-2)(s-1)$, with equality
if $I(X)$ is a complete intersection associated to a clutter with $s$
edges. 
\end{proposition}

\begin{remark} In Theorem~\ref{march14-10}(c) we provide another upper
bound for the index of regularity of $S/I(X)$ valid for a certain
family of algebraic toric sets. 
\end{remark}

\section{The degree-complexity of $I(X)$}\label{d-compl-section}
We continue to use the notation and definitions used in the
introduction. The main result of this section is an upper bound for
the degree-complexity of $I(X)$.  

In what follows we shall assume that $\succ$ is the {\it reverse
lexicographical order\/} ({\it revlex order\/} for short) on the
monomials of $S$. This order is given by
$t^b\succ t^a$ if and only if the last non-zero
entry  of $b-a$ is negative. As usual, if $g$ is a
polynomial of $S$, we denote the leading term of $g$ 
by ${\rm in}(g)$ and the leading coefficient of $g$ by ${\rm lc}(g)$. 

According to \cite[Proposition~6, p.~91]{CLO} the ideal $I(X)$ 
has a unique reduced Gr\"obner basis. We refer to \cite{CLO} for the
theory of Gr\"obner bases.  
The {\it degree-complexity\/} of $I(X)$, with respect to $\succ$, is
the maximum degree of the polynomials in the reduced 
Gr\"obner basis of $I(X)$. Next we study the reduced Gr\"obner basis
and the degree-complexity of $I(X)$. 

We come to one of the main results of this section. 

\begin{theorem}\label{gbasis-parameterized} Let $\mathcal{C}$ be a
clutter and let $\succ$ be the revlex
order on the monomials of $S$. If $\mathcal{G}$ is the reduced
Gr\"obner basis of the 
ideal $I(X)$, then $t_i^{q-1}-t_s^{q-1}\in \mathcal{G}$ for
$i=1,\ldots,s-1$ and $\deg_{t_i}(g)\leq q-1$ for $g\in\mathcal{G}$ 
and $1\leq i\leq s$.
\end{theorem}

\begin{proof} The reduced Gr\"obner basis of $I(X)$ consists of
homogeneous binomials \cite{algcodes}. 
As $I(X)$ is a lattice ideal \cite{algcodes}, it is seen
that each binomial $t^a-t^b\in\mathcal{G}$ satisfies that ${\rm
supp}(a)\cap{\rm supp}(b)=\emptyset$, this follows using that each 
variable $t_i$ is not a zero-divisor of $S/I(X)$. Since
$t_i^{q-1}-t_s^{q-1}$ is in $I(X)$ for 
$i=1,\ldots,s-1$, there is $g_i\in\mathcal{G}$ such that
$g_i=t_i^{b_i}-t^{c_i}$, $b_i\leq q-1$, $c_i\in\mathbb{N}^s$,
$i\notin{\rm supp}(c_i)$, and ${\rm in}(g_i)=t_i^{b_i}$. Then, by
Lemma~\ref{jan6-10-1}, the binomial $g_i$ has the form
$g_i=t_i^{q-1}-t_{j_i}^{q-1}$ for some $i<j_i$. As $\mathcal{G}$ is a
reduced Gr\"obner basis we get that $g_i=t_i^{q-1}-t_s^{q-1}$ for
$i=1,\ldots,s-1$. Let
$g\in\mathcal{G}\setminus\{g_1,\ldots,g_{s-1}\}$. Using that
$\mathcal{G}$ is reduced we get
that $\deg_{t_i}(g)\leq q-2$ for $i=1,\ldots,s-1$. To complete the
proof we need only show $\deg_{t_s}(g)\leq q-1$. Assume that
$a_s=\deg_{t_s}(g)>q-1$. After permuting $t_1,\ldots,t_{s-1}$ we may
assume that ${\rm in}(g)=t_1^{a_1}\cdots t_r^{a_r}$ and 
$g=t_1^{a_1}\cdots t_r^{a_r}-t_{r+1}^{a_{r+1}}\cdots t_s^{a_s}$, 
where $r<s$. Consider the polynomial
\begin{eqnarray*}
h&=& t_2^{a_2}\cdots t_r^{a_r}g_1-
t_1^{q-1-a_1}g\\
&=&t_s^{q-1}\left(-t_2^{a_2}\cdots
t_r^{a_r}+t_1^{q-1-a_1}t_{r+1}^{a_{r+1}}\cdots t_{s-1}^{a_{s-1}}
t_s^{a_s-(q-1)}\right)=t_s^{q-1}h_1.
\end{eqnarray*}
Since $h\in I(X)$ and using that $I(X)$ is a lattice ideal, we get
that the binomial 
$$h_1=-t_2^{a_2}\cdots
t_r^{a_r}+t_1^{q-1-a_1}t_{r+1}^{a_{r+1}}\cdots t_{s-1}^{a_{s-1}} 
t_s^{a_s-(q-1)}$$
belongs to $I(X)$. As ${\rm in}(h_1)=t_2^{a_2}\cdots
t_r^{a_r}$, we obtain that ${\rm in}(g)\in({\rm
in}(\mathcal{G}\setminus\{g\})$, a contradiction. Thus
$\deg_{t_s}(g)\leq q-1$. 
\end{proof}

The next result is interesting because it shows that the Hilbert
functions of $S/I(X)$ and $S/I_\mathcal{A}$ are equal up to degree
$q-2$. 

\begin{proposition} Let $\mathcal{C}$ be a clutter. If
$f=t^{a^+}-t^{a^-}$ is a non-zero homogeneous binomial of $I(X)$ 
and $\deg(f)\leq q-2$, then $f\in I_\mathcal{A}$.
\end{proposition}

\begin{proof} We may assume that $a^+=(a_1,\ldots,a_r,0,\ldots,0)$ 
and $a^{-}=(0,\ldots,0,a_{r+1},\ldots,a_m,0,\ldots,0)$ and $a_i\geq 1$
for $i=1,\ldots,m$. Then 
\begin{equation}\label{jan12-10}
(x_1^{v_{11}}\cdots
x_n^{v_{1n}})^{a_1}\cdots (x_1^{v_{r1}}\cdots
x_n^{v_{rn}})^{a_r}=(x_1^{v_{r+1,1}}\cdots
x_n^{v_{r+1,n}})^{a_{r+1}}\cdots(x_1^{v_{m,1}}\cdots
x_n^{v_{m,n}})^{a_{m}}
\end{equation}
for all $(x_1,\ldots,x_n)\in (K^*)^n$, where 
$v_i=(v_{i1},\ldots,v_{in})=(v_{i,1},\ldots,v_{i,n})$. To show that
$f\in I_\mathcal{A}$ we need only show that $A{a^+}=A{a^-}$, where $A$
is the incidence matrix of $\mathcal{C}$, i.e., $A$ is the matrix with column
vectors $v_1,\ldots,v_s$. Equivalently we need only show the equality 
\begin{equation}\label{jan12-10-1}
v_{1,k}a_1+\cdots+v_{r,k}a_r=v_{r+1,k}a_{r+1}+\cdots+v_{m,k}a_m
\end{equation}
for $1\leq k\leq n$. If both sides of Eq.~(\ref{jan12-10-1}) are zero
there is nothing to show. We proceed by contradiction assuming:
\begin{equation}\label{jan12-10-2}
v_{1,k}a_1+\cdots+v_{r,k}a_r>v_{r+1,k}a_{r+1}+\cdots+v_{m,k}a_m\geq 0.
\end{equation}
Making $x_i=1$ for $i\neq k$ in Eq.~(\ref{jan12-10}), we get 
$$
x_k^{v_{1,k}a_1+\cdots+v_{r,k}a_r}=x_k^{v_{r+1,k}a_{r+1}+\cdots+v_{m,k}a_m}
$$
for any $x_k\in K^*$. In particular making $x_k=\beta$, where $\beta$
is a generator of the cyclic group $(K^*,\,\cdot\, )$, we get 
that
\begin{equation}\label{jan9-11}
(v_{1,k}a_1+\cdots+v_{r,k}a_r)-(v_{r+1,k}a_{r+1}+\cdots+v_{m,k}a_m)
\equiv\ 0\ {\rm mod} (q-1). 
\end{equation}
Consequently $v_{1,k}a_1+\cdots+v_{r,k}a_r\geq q-1$, a contradiction
because 
$$
q-2\geq \deg(f)=a_1+\cdots+a_r\geq v_{1,k}a_1+\cdots+v_{r,k}a_r.
$$ 
Hence equality in Eq.~(\ref{jan12-10-1}) holds for $1\leq k\leq n$ and
the proof is complete.
\end{proof}

\begin{proposition} Let $A$ be the matrix with column vectors
$v_1,\ldots,v_s$. Then
$$I(X)=(\{t^{a^+}-t^{a^-}\vert\, Aa^+ \equiv Aa^-\, {\rm mod}\, (q-1)\mbox{
and }|a^+|=|a^-|\}).$$
\end{proposition}

\begin{proof} The inclusion ``$\subset$'' follows from
Eq.~(\ref{jan9-11}) and from the fact that $I(X)$ is a lattice ideal
\cite{algcodes}. To show the inclusion ``$\supset$'' take
$f=t^{a^+}-t^{a^-}$ such that $Aa^+ \equiv Aa^-\, {\rm mod}\, (q-1)$
and $|a^+|=|a^-|$. From the first condition it is seen that $f$
vanishes on $X$ and from the second condition $f$ is homogeneous
in the standard grading of $S$. Thus $f\in I(X)$.
\end{proof}

\section{Upper bounds for the minimum
distance}\label{upper-bounds-section}

We continue to use the notation and definitions used in the
introduction and in the preliminaries. Let $\mathcal{C}$ be a clutter
with vertex set 
$V_\mathcal{C}=\{y_1,\ldots,y_n\}$. Throughout
this section we assume that $\mathcal{A}=\{v_1,\ldots,v_s\}$ is the
set of all 
characteristic vectors of the edges of a uniform clutter $\mathcal{C}$.

The set $(K^*)^n$ is called an {\it
affine algebraic torus} of dimension $n$ and is denoted  by
$\mathbb{T}^*$. The torus $\mathbb{T}^*$ is a multiplicative group under
the product operation $(\alpha_i)(\alpha_i')=(\alpha_i\alpha_i')$, where
$(\alpha_i)$ really means $(\alpha_1,\ldots,\alpha_n)$. Clearly, the
algebraic toric set: 
$$
X:=\{[(x_1^{v_{11}}\cdots x_n^{v_{1n}},\ldots,x_1^{v_{s1}}\cdots
x_n^{v_{sn}})]\vert\, x_i\in K^*\mbox{ for all
}i\}\subset\mathbb{P}^{s-1}
$$
is also a multiplicative group with the product operation. 

Let $I$ be the ideal of $R=K[y_1,\ldots,y_n]$ generated by
$y^{v_1},\ldots,y^{v_s}$. The ideal $I$ is called the {\it edge ideal} of
$\mathcal{C}$ and the matrix $A$ whose columns are $v_1,\ldots,v_s$ is
called the {\it incidence matrix} of $\mathcal{C}$. Recall that the
{\it integral  
closure\/} of $I^i$, denoted by
$\overline{I^i}$, is the ideal of $R$ given by  
\begin{equation}\label{oct12-08}
\overline{I^i}=(\{y^a\in R\vert\, \exists\,
p\in\mathbb{N}\setminus\{0\};(y^a)^{p}\in
I^{pi}\}),
\end{equation}
see for instance \cite[Proposition~7.3.3]{monalg}. The ideal $I$ 
is called {\it normal\/} if $I^i=\overline{I^i}$ for $i\geq 1$. There
are many interesting examples of normal ideals \cite{Stur1,monalg}.
For instance if $\mathcal{C}$ is the clutter of all subsets of 
$Y=\{y_1,\ldots,y_n\}$ of a fixed size $k\geq 1$, then $I$ is normal.
If $\mathcal{C}$ is 
the clutter of bases of a matroid, then $I$ is also normal
\cite{matrof}. There is a
combinatorial description of the normality of ideals generated by
square-free monomials of degree $2$ \cite{bowtie}, i.e., of ideals such that
$\mathcal{C}$ is a graph. According to
\cite{bowtie} if $\mathcal{C}$ is a complete graph or a bipartite graph, then $I$ is normal.
The edge ideal $I$ is also normal if $\mathcal{C}$ is any odd cycle or
any unicyclic graph. 

Let $\mathcal{B}\subset\mathbb{Z}^{n+1}$. The {\it polyhedral
cone\/} generated by $\mathcal{B}$ is denoted by
$\mathbb{R}_+\mathcal{B}$. A polyhedral cone containing no lines is called 
{\it pointed}. The set $\mathcal{B}$ is 
called a {\it Hilbert basis} if
$\mathbb{N}\mathcal{B}=\mathbb{R}_+\mathcal{B}\cap\mathbb{Z}^{n+1}$,
where $\mathbb{N}\mathcal{B}$ is the semigroup generated by
$\mathcal{B}$. 

We come to the main result of this section, an upper bound for the
minimum distance of $C_X(d)$ valid for certain normal edge ideals of
uniform 
clutters. 

\begin{theorem}\label{march14-10} Let $\mathcal{C}$ be a uniform
clutter whose incidence matrix 
has rank $n$ and let 
$I\subset R$ be its edge ideal. If $I$ is normal and $\mathbb{T}$ is
a projective  
torus in $\mathbb{P}^{n-1}$, then:
\begin{enumerate} 
\item[\rm (a)] The degree of $S/I(X)$ is equal to $|X|=(q-1)^{n-1}$.
\item[\rm (b)] $\delta_d\leq\delta_d'$, where $\delta_d'$ is the
minimum distance of the linear code $C_\mathbb{T}(d)$.
\item[\rm (c)] ${\rm reg}(S/I(X))\leq {\rm
reg}(S'/I(\mathbb{T}))=(q-2)(n-1)$,
where $S'=K[t_1,\ldots,t_n]$.
\end{enumerate}
\end{theorem}

\begin{proof} (a) The ideal $I$ is normal. Then by \cite[Theorem~3.15]{ehrhart} the set
$\mathcal{B}=\{(v_i,1)\}_{i=1}^s$ is a Hilbert basis. Therefore, using
\cite[Theorem~3.5]{algcodes}, we obtain that $(q-1)^{n-1}$ divides
$|X|$. On the other hand there is an epimorphism of multiplicative groups 
$$
\theta\colon\mathbb{T}^*\rightarrow X\, ;\ \ \ \ \ \ \ \ 
(x_1,\ldots,x_n)\stackrel{\theta}{\longmapsto}
[(x^{v_1},\ldots,x^{v_s})],
$$ 
where $\mathbb{T}^*=(K^*)^{n}$ is an affine algebraic torus. The kernel of $\theta$
contains the diagonal  
subgroup 
$$\mathcal{D}^*=\{(\lambda,\ldots,\lambda)\vert\, \lambda\in
K^*\}.$$
Thus $|X|$ divides $(q-1)^{n-1}$. Putting altogether, we get
$|X|=(q-1)^{n-1}$. 

(b) The set $\mathcal{B}=\{(v_i,1)\}_{i=1}^s$
is a Hilbert basis (see the proof of part (a)). Hence using a result of
\cite{gerards-sebo}, 
after permutation of the $(v_i,1)$'s,  we may assume that 
$\mathcal{B}'=\{(v_1,1),\ldots,(v_n,1)\}$ is a Hilbert basis and a
linearly independent set. Then, it is seen that the group
$\mathbb{Z}^{n+1}/\mathbb{Z}\mathcal{B}'$ 
is torsion-free. We set $L'=\mathbb{Z}\{v_i-v_1\}_{i=2}^n$. 
It is
not hard to see that there is an isomorphism of groups 
$$
\tau{\colon}T(\mathbb{Z}^n/L')\rightarrow 
T(\mathbb{Z}^{n+1}/\mathbb{Z}\mathcal{B}')
$$ 
given by $\tau(\overline{a})=\overline{(a,0)}$, where $T(M)$
denotes the torsion subgroup of an abelian group $M$, i.e., $T(M)$, is the set of
all $m$ in $M$ such that $pm=0$
for some $0\neq p\in \mathbb{Z}$. From this isomorphism we conclude
that $T(\mathbb{Z}^n/L')=0$, i.e., $\mathbb{Z}^n/L'$ is also
torsion-free. 

Consider the algebraic toric set parameterized by $y^{v_1},\ldots,
y^{v_n}$:
$$X_1=\{[(x^{v_1},\ldots, x^{v_n})]\vert\, x_i\in
K^*\mbox{ for all
}i\}\subset\mathbb{P}^{n-1}.
$$ 
We claim that $I(X_1)=(\{t_i^{q-1}-t_n^{q-1}\}_{i=1}^{n-1})$. We set
$\mathcal{A}'=\{v_1,\ldots,v_n\}$. Notice that the set $\mathcal{A}'$
is also linearly independent. 
Since $I_{\mathcal{A}'}=(0)$ and
$\mathbb{Z}^{n}/L'$ is torsion-free, 
by Theorem~\ref{vila-dictaminadora}, we obtain
\begin{eqnarray*}
(\{t_i^{q-1}-t_n^{q-1}\}_{i=1}^{n-1})&=&(\{t_i^{q-1}-t_n^{q-1}\}_{i=1}^{n-1})\colon (t_1\cdots
t_n)^{\infty})\\
&=&(I_{\mathcal{A}'}+(\{t_i^{q-1}-t_n^{q-1}\}_{i=1}^{n-1}))\colon (t_1\cdots
t_n)^{\infty})\\ 
&\stackrel{\ref{vila-dictaminadora}}{=}&I(X_1).
\end{eqnarray*}

Let $\mathbb{T}$ be a projective torus in 
$\mathbb{P}^{n-1}$. By Proposition~\ref{ci-summary}, we have
$I(\mathbb{T})=I(X_1)$. Consequently $X_1=\mathbb{T}$
because $X_1$ and $\mathbb{T}$ are projective varieties. Let $\delta_d'$ be the
minimum distance of $C_{X_1}(d)$. Next we show that
$\delta_d\leq\delta_d'$. There is a
well defined epimorphism $$
\overline{\theta}_1\colon X\rightarrow X_1,\ \ \ \ \ \ \ \
[(x^{v_1},\ldots,x^{v_s})]\mapsto [(x^{v_1},\ldots,x^{v_n})] 
$$
induced by the projection map
$[(\alpha_1,\ldots,\alpha_s)]\mapsto[(\alpha_1,\ldots,\alpha_n)]$. By part (a) one has
$|X|=|X_1|=(q-1)^{n-1}$. Therefore the map $\overline{\theta}_1$ 
is an isomorphism of multiplicative groups. For any homogeneous
polynomial $F$, we denote its zero set by $Z_X(F)=\{[P]\in X\, \vert\, F(P)=0\}$.
Let $S'=K[t_1,\ldots,t_n]=\oplus_{d=0}^\infty S_d'$ and let $F_1\in
S_d'$ be a polynomial such that ${\rm ev}_d(F_1)\neq 0$ and with
$|Z_{X_1}|$ as large as possible, i.e., we choose $F_1$ so that
$\delta_d'=|X_1|-|Z_{X_1}(F_1)|$. We can regard the polynomial
$F_1=F_1(t_1,\ldots,t_n)$ as an element of $S$ and denote it by $F$. 
The map $\overline{\theta}_1$
induces a bijective map
$$
\overline{\theta}_1\colon Z_X(F)\mapsto Z_{X_1}(F_1),\ \ \ \ \ \ \ \ 
[P]\mapsto [\overline{\theta}_1(P)].
$$     
Therefore we have the inequality
\begin{eqnarray*}
\max\{|Z_X(F)|\colon F\in
S_d;\, {\rm ev}_d(F)\neq 0\}&\geq &\max\{|Z_{X_1}(F_1)|\colon F_1\in
S_d';\, {\rm ev}_d(F_1)\neq 0\}.
\end{eqnarray*}
Consequently  $\delta_d\leq \delta_d'$. 

(c) We continue to use the notation and definitions used in the proof
of part (b). Since $X_1=\mathbb{T}$, it suffices to show that
$H_{X_1}(d)\leq H_X(d)$ for $d\geq 1$. Using that $I(X)$ and $I(X_1)$
are vanishing ideals generated by homogeneous polynomials, it is not
hard to show that $S'\cap I(X)=I(X_1)$. Thus, we have a graded monomorphism
$$
0\rightarrow K[t_1,\ldots,t_n]/I(X_1)\rightarrow
K[t_1,\ldots,t_s]/I(X),\ \ \ \ \ \overline{F}_1\mapsto \overline{F}_1.
$$
Hence $H_{X_1}(d)\leq H_X(d)$ and consequently ${\rm reg}(S/I(X))\leq
{\rm reg}(S'/I(\mathbb{T}))=(q-2)(n-1)$.  
\end{proof}

There is a nice recent formula for $\delta_d'$:  

\begin{theorem}{\cite[Theorem~3.4]{ci-codes}}\label{maria-vila-hiram-eliseo} 
If $\mathbb{T}$ is a projective torus in $\mathbb{P}^{n-1}$ and $d\geq 1$, 
then the minimum distance of $C_\mathbb{T}(d)$ is given by 
$$
\delta_d'=\left\{\begin{array}{cll}
(q-1)^{n-(k+2)}(q-1-\ell)&\mbox{if}&d\leq (q-2)(n-1)-1,\\
1&\mbox{if}&d\geq (q-2)(n-1),
\end{array}
 \right.
$$
where $k$ and $\ell$ are the unique integers such that $k\geq 0$,
$1\leq \ell\leq q-2$ and $d=k(q-2)+\ell$. 
\end{theorem}

\begin{remark}\label{comparison-remark-1}\rm (i) When $d$ is greater
than or equal to the index of regularity of 
$S/I(X)$, by Proposition~\ref{request-referee}, one has that
$\delta_d=1$. Thus,  
for $d\geq {\rm reg}(S/I(X))$ our codes are useless from a practical point of
view. For some other values of the parameters
however, the bound $\delta_d'$ does not prevent our codes from having a large
(although not optimal) minimum distance. In Example~\ref{large-minimum-distance} we provide
specific values of the parameters of $C_X(d)$ when $\mathcal{C}$ is a cycle of
length $3$. 

(ii) Let $\mathcal{C}$ be a unicyclic connected graph with $n$
vertices and with a unique
cycle of odd length.  Then, $X=\mathbb{T}$ is a
projective torus in 
$\mathbb{P}^{n-1}$ by Corollary~\ref{jan4-11}. Thus, the minimum
distance $\delta_d$ of $C_X(d)$ is equal to 
$\delta_d'$ by Theorem~\ref{maria-vila-hiram-eliseo}. In particular
$\delta_d=1$ for $d\geq (q-2)(n-1)$.

(iii) The problem of computing the minimum distance of a linear code is
NP-hard \cite{vardy-complexity}. It might not be easy to 
compute the minimum distance of $C_X(d)$ for graphs with large number of edges and
vertices. However, for a complete graph with $4$ vertices it is not hard to
compute the minimum distance and to compare the bound $\delta_d'$
with 
the Singleton bound, see
Example~\ref{k4-compare-bounds}.   

\end{remark}

\begin{example}\label{large-minimum-distance}
Let $\mathcal{C}$ be a cycle of length $3$, let $X$ be the
algebraic toric set parameterized by $y_1y_2,y_2y_3,y_1y_3$ and
let $C_X(d)$ be the parameterized code of order $d$ over the field
$K=\mathbb{F}_{9}$. Using {\em Macaulay\/}$2$, together with
Remark~\ref{comparison-remark-1}(ii), we obtain the basic
parameters of $C_X(d)$:

\begin{eqnarray*}
&&\left.
\begin{array}{c|c|c|c|c|c|c|c|c|c|c|c|c|c|c}
 d & 1 & 2 & 3 & 4 & 5 & 6 & 7 & 8 & 9 & 10 & 11 & 12 & 13 & 14  \\
   \hline
 |X| & 64 & 64 & 64 & 64 & 64 & 64 & 64 & 64 & 64 & 64 & 64 & 64 & 64
 & 64  \\ 
   \hline
 \dim C_X(d)    \    & 3 & 6   & 10 & 15 & 21 & 28 & 36 & 43 & 49 &
 54 & 58 & 61 & 63 & 64  \\ 
   \hline
 \delta_d & 56 & 48 & 40 & 32 & 24 & 16 & 8 & 7 & 6 & 5  & 4  & 3   &
 2  & 1 \\ 
\end{array}
\right.
\end{eqnarray*}

For linear codes over $\mathbb{F}_q$ with $q\leq 9$, there are online
tables of known upper and lower bounds on the optimal minimum
distance for each given dimension \cite{codetables}. The last
line of the following table shows the upper bounds for the minimum
distance of $C_X(d)$ that we found using \cite{codetables}.

\begin{eqnarray*}
&&\left.
\begin{array}{c|c|c|c|c|c|c|c|c|c|c|c|c|c|c}
 d & 1 & 2 & 3 & 4 & 5 & 6 & 7 & 8 & 9 & 10 & 11 & 12 & 13 & 14  \\
   \hline
 |X| & 64 & 64 & 64 & 64 & 64 & 64 & 64 & 64 & 64 & 64 & 64 & 64 & 64
 & 64  \\ 
   \hline
 \dim C_X(d)    \    & 3 & 6   & 10 & 15 & 21 & 28 & 36 & 43 & 49 &
 54 & 58 & 61 & 63 & 64  \\ 
   \hline
 m_d & 56 & 53 & 49 & 44 & 39 &  32& 26 & 19 & 13 & 9  & 5  & 3   &
 2  & 1 \\ 
\end{array}
\right.
\end{eqnarray*}

The $C_X(d)$ linear codes for this example are really
only competitive---with other known codes of the same block
length and dimension---in the very low rate cases (i.e. small $d$
where the dimension is much less 
than the length) and the very high rate cases (i.e. $d$ close to
$(q-2)(n-1)$). 
\end{example}

\begin{example}\label{k4-compare-bounds}
Let $\mathcal{C}$ be the following complete graph on four vertices and let
$X$ be the algebraic toric set parameterized by all $y_iy_j$ such
that $\{y_i,y_j\}$ is an edge of $\mathcal{C}$.

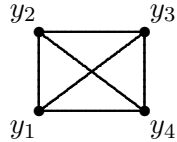
\bigskip

\begin{picture}(230,50)(-100,0)
\thicklines
 \put(90,10){\circle*{4}}
 \put(84,3){\makebox(0,0)[c]{$y_1$}}
 \put(90,40){\circle*{4}}
 \put(84,47){\makebox(0,0)[c]{$y_2$}}
 \put(130,40){\circle*{4}}
 \put(137,47){\makebox(0,0)[c]{$y_3$}}
 \put(130,10){\circle*{4}}
 \put(137,3){\makebox(0,0)[c]{$y_{4}$}}
 \curve(90,10, 90,40)
 \curve(90,10, 130,40)
 \curve(90,10, 130,10)
 \curve(90,40, 130,40)
 \curve(90,40, 130,10)
 \curve(130,40, 130,10)
\end{picture}

\medskip

Let $C_X(d)$ be the parameterized code of order $d$ over the field
$K=\mathbb{F}_3$, let $b_d$ (resp. $\delta_d'$) be the Singleton
bound (resp. the bound of Theorem~\ref{march14-10}), and let
$\delta_d$ be the minimum distance
of $C_X(d)$. Using {\em Macaulay\/}$2$, we obtain:
\begin{eqnarray*}
&&\left.
\begin{array}{c|c|c|c}
 d & 1 & 2 & 3  \\ \hline
 b_d & 3 & 1 & 1 \\ \hline
 \delta_d'& 4 & 2 & 1  \\ \hline
 \delta_d & 2 & 1 & 1 
\end{array}\right.
\end{eqnarray*}

\medskip

If $C_X(d)$ is the parameterized code of order $d$ over the
field $K=\mathbb{F}_4$, then we get:
\begin{eqnarray*}
&&\left.
\begin{array}{c|c|c|c|c|c|c}
 d & 1 & 2 & 3 & 4 & 5 & 6\\ \hline
 b_d & 22 & 9 & 1& 1 & 1& 1 \\ \hline 
 \delta_d'&  18 & 9 & 6 & 3 & 2 & 1 \\ \hline
 \delta_d & 12 & 3 & 1 & 1 & 1 & 1
\end{array}\right.
\end{eqnarray*}
\end{example}

\medskip

\begin{center}
ACKNOWLEDGMENT
\end{center}

\noindent We thank Hiram L\'opez and Carlos Renter\'\i a 
for many stimulating discussions. We also thank the referees for their
careful reading of the paper and for the improvements that
they suggested.

\bigskip

\bibliographystyle{plain}

\end{document}